\documentclass[12pt]{amsart}
\usepackage{amsmath, amssymb, amsthm}
\usepackage{paralist, xcolor, hyperref}

\usepackage[margin=1in,letterpaper,portrait]{geometry}
\usepackage{diagbox}
\usepackage[capitalize]{cleveref}

\usepackage{tikz}
\usetikzlibrary{patterns}

\newcommand{\inv}{I}
\newcommand{\invol}{\mathcal{I}}

\newcommand{\Int}{\mathrm{Int}}
\newcommand{\val}{\nu}
\newcommand{\poset}{\mathcal{P}}
\newcommand{\0}{\hat{0}}
\newcommand{\1}{\hat{1}}
\newcommand{\id}{e_n}

\newtheorem{conjecture}{Conjecture}[section]
\newtheorem{theorem}[conjecture]{Theorem}
\newtheorem{corollary}[conjecture]{Corollary}
\newtheorem{lemma}[conjecture]{Lemma}
\newtheorem{proposition}[conjecture]{Proposition}

\theoremstyle{definition}

\newtheorem{example}[conjecture]{Example}
\newtheorem{remark}[conjecture]{Remark}
\newtheorem{definition}[conjecture]{Definition}
\newtheorem*{question*}{Question}

\usepackage{todonotes}

\title[The middle order on permutations]{Between weak and Bruhat: \\the middle order on permutations}

\author{Mathilde Bouvel}
\address{Université de Lorraine, CNRS, Inria, LORIA, F-54000 Nancy, France}
\email{mathilde.bouvel@loria.fr}

\author[Luca Ferrari]{Luca Ferrari$^\ddagger$}
\address{Department of Mathematics and Computer Science, University of Firenze, Firenze, Italy}
\email{luca.ferrari@unifi.it}
\thanks{$^\ddagger$Research partially supported by INdAM -- GNCS project CUP\_E53C23001670001.}

\author[Bridget Eileen Tenner]{Bridget Eileen Tenner$^\dagger$}
\address{Department of Mathematical Sciences, DePaul University, Chicago, IL, USA}
\email{bridget@math.depaul.edu}
\thanks{$^\dagger$Research partially supported by NSF Grant DMS-2054436.}

\begin{document}

\begin{abstract}
We define a partial order $\poset_n$ on permutations of any given size $n$, which is the image of a natural partial order on inversion sequences. 
We call this the ``middle order.'' We demonstrate that the poset $\poset_n$ refines the weak order on permutations and admits the Bruhat order as a refinement, justifying the terminology. These middle orders are distributive lattices and we establish some of their combinatorial properties, including characterization and enumeration of intervals and boolean intervals (in general, or of any given rank), and a combinatorial interpretation of their Euler characteristic. We further study the (not so well-behaved) restriction of this poset to involutions, obtaining a simple formula for the M\"obius function of  principal order ideals there. Finally, we offer further directions of research, initiating the study of the canonical Heyting algebra associated with $\poset_n$, and defining a parking function analogue of~$\poset_n$. 
\end{abstract}

\maketitle

In this paper, we are interested in partial orders that can be defined on the sets of permutations of any fixed size. 
Among the most classical such orders are the Bruhat order and the weak order, which motivate the context of our study. 
These are defined on the set $S_n$ of all permutations of size $n$, and can be conveniently described in terms of the associated covering relations. More specifically, given $v,w\in S_n$, the element $v$ is covered by $w$ in the Bruhat order when $w$ is obtained from $v$ by turning a noninversion into an inversion, provided that no further inversions are created. 
On the other hand, $v$ is covered by $w$ in the weak order when $w$ is obtained from $v$ by turning an ascent into a descent. 
We refer to~\cite{BB} for background on the Bruhat and weak orders, and to \cite{ec1} for basics on posets in general. 

Among the interesting properties of these orders, we remark that the weak orders are lattices, whereas the Bruhat orders are not lattices in general. For a generic poset, the lattice property may fail because there are pairs of elements that either have no upper or lower bounds, or do have those bounds but do not have a \emph{least} upper or \emph{greatest} lower bound. Since the Bruhat order on $S_n$ clearly has a minimum (the identity permutation $\id$) and a maximum (the reverse identity permutation, also called the \emph{long element}), the lattice property fails in this poset because of the second reason listed above, and this is somehow related to the fact that there are too many covering relations. Because the Bruhat orders are refinements of the weak orders, we may assert that, in some sense, the weak orders are lattices because enough covering relations have been removed from the Bruhat orders. Note, as will become relevant soon, that the weak orders are not generally distributive as lattices. 

In the present paper, we introduce a class of ``intermediate'' partial orders, which interpolate (with respect to refinement) between the Bruhat orders and the weak orders. Among other features, these new posets, which we call \emph{middle orders}, 
have the appealing property of being distributive lattices. 
In addition, while the Bruhat and weak orders can be defined in terms of the reduced decompositions of a permutation, the middle order is naturally defined in terms of inversion sequences associated with permutations; to our knowledge, this is the first time that a connection between inversion sequences and partial orders on permutations is studied. 
Our efforts here are to study fundamental combinatorial properties of the middle order on permutations, and to establish its combinatorial relevance. The details of this work lead us to suspect that properties of this poset may well have uses in and applications to other fields.

In Section~\ref{middle} we introduce the middle order on permutations, and discuss its relationship with the weak and the Bruhat orders. In particular, we provide a characterization of the covering relations in terms of mesh patterns. Section~\ref{sec:ipf} links our work to previous literature, where the middle order has occasionally surfaced (under different names). Section~\ref{sec:counting} is devoted to the study of the enumerative combinatorics of intervals in general, and boolean intervals in particular. In Section~\ref{euler_char} we recall the notion of valuation in a distributive lattice, and we give a combinatorial description for the Euler characteristic of the middle order. Restricting to involutions gives rise to subposets that are much more difficult to investigate (in particular, they are not lattices in general). For such subposets, we address and completely solve the problem of computing the M\"obius function of  principal order ideals, which is the content of Section~\ref{sec:involutions}. We conclude with Section~\ref{open_problems}, containing suggestions for further research.

\section{The middle order}\label{middle}

An \emph{inversion} in a permutation is an out-of-order pair in the one-line notation of the permutation. 
For example, in the permutation $415623$, the third entry $5$ (sometimes called ``inversion top'') and the fifth entry $2$ (sometimes called ``inversion bottom'') form an inversion. 
Inversions are often described in terms of the indices of the entries in that out-of-order pair, but of course it is equivalent to describe their values instead. 
Accordingly, there are several ways to define the \emph{inversion sequence} of a permutation. Here, we use the following.

\begin{definition}\label{defn:inversion sequence}
For a permutation $w \in S_n$, the \emph{inversion sequence} of $w$ is
$$\inv(w) = (x_1, \ldots, x_n),$$
where
$$x_i = \#\{j< i \, | \, w^{-1}(j) > w^{-1}(i)\}.$$
\end{definition}

In other words, the $i$th coordinate of $\inv(w)$ counts the times that the value $i$ is an inversion top. 
Going back to our example of $w = 415623$, we have $I(w)=(0,0,0,3,2,2)$; for instance, $x_5 =2$ since the value $5$ is the inversion top of two inversions (whose inversion bottoms are $2$ and $3$). 
The main object of our work is a new partial order on $S_n$, which can be expressed naturally in terms of these inversion sequences, as presented in the following definition.

\begin{definition}\label{defn:poset}
    Let $\poset_n$ be the poset whose elements are the permutations in $S_n$, with ordering relations defined by $v \le w$ if and only if the inversion sequence $\inv(v)$ is, coordinate-wise, less than or equal to the inversion sequence $\inv(w)$. We call $\poset_n$ the \emph{middle order} of size $n$, and we denote the coordinate-wise-ordering of inversion sequences by $I(v) \le I(w)$. 
\end{definition}

The covering relations in $\poset_n$ are characterized by 
\begin{eqnarray*}
   & v \lessdot w &\\
   & \Updownarrow &\\
   & \inv(w) = \inv(v) + (0,\ldots, 0,1,0,\ldots,0) &
\end{eqnarray*} 
where the second term in the sum is a standard basis vector of appropriate size $n$.

Clearly an inversion sequence $(x_1, \ldots, x_n)$ as in Definition~\ref{defn:inversion sequence} must satisfy $x_i \in [0, i-1]$ for each $i$. Moreover, it is well known that $S_n$ is in bijection with the set of sequences $\{(y_1, \ldots, y_n)\, | \, y_i \in [0,i-1]\}$, via the correspondence $w \mapsto \inv(w)$.
We will let $I^{-1} : [0,0] \times [0,1] \times [0,2] \times \cdots \times [0,n-1] \rightarrow S_n$ denote the map that is the inverse of this correspondence. 
This gives a natural isomorphism between posets.

\begin{proposition}\label{prop:poset bijection}
    The middle order poset $\poset_n$ is isomorphic to the product of chains
    $$[0,0] \times [0,1] \times [0,2] \times \cdots \times [0,n-1].$$
\end{proposition}
\begin{proof}
    The bijective correspondence between $\poset_n$ and $[0,0] \times [0,1] \times [0,2] \times \cdots \times [0,n-1]$ has been described above. 
    The poset isomorphism follows since, by definition of the middle order, $v \leq w$ if and only if $I(v)$ is less than or equal to $I(w)$ coordinate-wise. 
\end{proof}

A poset class of particular interest is the collection of lattices.
More information on these objects can be found, for example, in \cite[Chapter 3]{ec1}.

\begin{definition}\label{defn:lattice meet join}
    A poset is a \emph{lattice} if every pair of elements $s$ and $t$ has a least upper bound, denoted $s \vee t$, and a greatest lower bound, denoted $s \wedge t$. The former of these is called the \emph{join} of $s$ and $t$, and the latter is the \emph{meet}.
\end{definition}

The poset $\poset_n$ is a lattice, and the meet and join operations are characterized  coordinate-wise on the corresponding inversion sequences. More precisely, for $v,w \in \poset_n$ such that $\inv(v) = (x_1,\ldots, x_n)$ and $\inv(w) = (y_1,\ldots, y_n)$, we have 
\begin{align}
    v \vee w &= I^{-1}((\max\{x_1,y_1\},\ldots, \max\{x_n,y_n\})) \text{ and} \label{eq:meet_and_join1} \\
    v \wedge w &= I^{-1}((\min\{x_1,y_1\},\ldots, \min\{x_n,y_n\})) \label{eq:meet_and_join2}, 
\end{align}
where $I^{-1}$ is the map from inversion sequences to permutations defined above. 
We also note that the minimum element in $\poset_n$ is the identity permutation $\id$, whose inversion sequence is the zero vector. 

\begin{definition}\label{defn:distributive}
    A lattice is a \emph{distributive} lattice if the meet and join operations distribute over each other. That is, if for all $s$, $t$, and $u$ in the poset, we have
    \begin{align*}
        s \vee (t \wedge u) &= (s \vee t) \wedge (s \vee u) \text{ \ and}\\
        s \wedge (t \vee u) &= (s \wedge t) \vee (s \wedge u).
    \end{align*}
\end{definition}

It follows from Proposition~\ref{prop:poset bijection} that $\poset_n$ is a finite distributive lattice. One consequence of that property is that a ``rank function'' can be defined on $\poset_n$.

\begin{definition}\label{defn:rank function}
    A poset is \emph{graded} if every maximal chain in the poset has the same length. Equivalently, a poset is graded if there exists a function $r$ on the elements of the poset where $r(x)$ is equal to the length of any maximal chain from $x$ to a minimal element of the poset. This $r$ is the \emph{rank function} of the poset. 
\end{definition}

The rank function is defined on the elements of a graded poset, but can be extended to intervals of this poset; specifically, for two elements $x,y$ with $x \leq y$ in the poset, the rank of the interval $[x,y]$ is $r([x,y]) := r(y) -r(x)$.

Finite distributive lattices are graded, so $\poset_n$ has a rank function. The rank $r(w)$ of any element $w$ in $\poset_n$ has a simple combinatorial interpretation as the number of inversions of $w$.

\begin{figure}[htbp]
\begin{center}
\begin{tikzpicture}
    \draw (0,0) coordinate (123);
    \draw (1.5,1) coordinate (132);
    \draw (-1.5,1) coordinate (213);
    \draw (1.5,2.5) coordinate (312);
    \draw (-1.5,2.5) coordinate (231);
    \draw (0,3.5) coordinate (321);
    \foreach \x in {213,132} {\draw (123) -- (\x); \draw (231) -- (\x); \draw (312) -- (\x);}
    \foreach \x in {231, 312} {\draw (321) -- (\x);}
    \foreach \x in {123,132,213,231,312,321} {\fill[white] (\x)++(-.35,-.25) rectangle ++(.7,.5); \draw (\x) node {$\x$};}
\end{tikzpicture}
    \hspace{.5in}
\begin{tikzpicture}
    \draw (0,0) coordinate (123);
    \draw (1.5,1) coordinate (132);
    \draw (-1.5,1) coordinate (213);
    \draw (1.5,2.5) coordinate (312);
    \draw (-1.5,2.5) coordinate (231);
    \draw (0,3.5) coordinate (321);
    \foreach \x in {213,132} {\draw (123) -- (\x); \draw (231) -- (\x);}
    \draw (132) -- (312);
    \foreach \x in {231, 312} {\draw (321) -- (\x);}
    \foreach \x in {123,132,213,231,312,321} {\fill[white] (\x)++(-.35,-.25) rectangle ++(.7,.5); \draw (\x) node {$\x$};}
\end{tikzpicture}
    \hspace{.5in}
\begin{tikzpicture}
    \draw (0,0) coordinate (123);
    \draw (1.5,1) coordinate (132);
    \draw (-1.5,1) coordinate (213);
    \draw (1.5,2.5) coordinate (312);
    \draw (-1.5,2.5) coordinate (231);
    \draw (0,3.5) coordinate (321);
    \foreach \x in {213,132} {\draw (123) -- (\x);}
    \draw (132) -- (312);
    \draw (213) -- (231);
    \foreach \x in {231, 312} {\draw (321) -- (\x);}
    \foreach \x in {123,132,213,231,312,321} {\fill[white] (\x)++(-.35,-.25) rectangle ++(.7,.5); \draw (\x) node {$\x$};}
\end{tikzpicture}
\end{center}
\caption{From left to right: the Bruhat order, the middle order and the weak order on $S_3$.}\label{biw}
\end{figure}
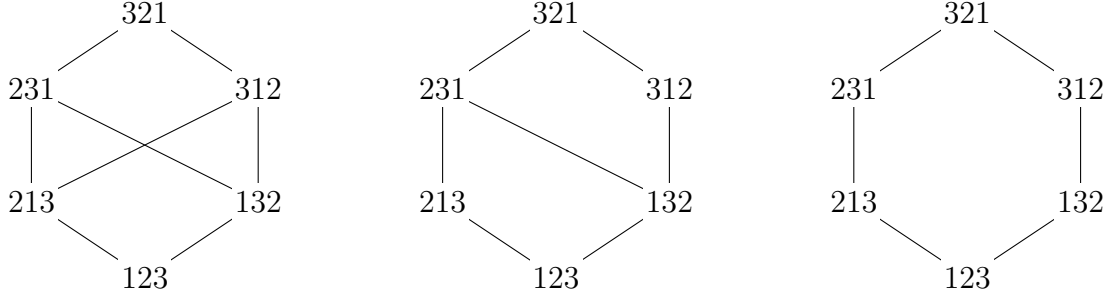

Figure~\ref{biw} illustrates the Hasse diagrams of the Bruhat, middle, and weak orders on $S_3$. As can be immediately noticed, each order is a refinement of the orders on its right. This is not an accident, as we will now prove. To this aim, we need to introduce some definitions and notations related to permutation patterns.

There is extensive literature related to patterns in permutations -- see the surveys~\cite{BonaBook,LivreKitaev,VatterSurvey}. 
The so-called ``mesh patterns,'' introduced by Br\"{a}nd\'en and Claesson in \cite{branden claesson}, are relevant to our work here.

\begin{definition}
    Fix a permutation $p \in S_k$ and a subset $M$ of unit squares in the grid $[0,k+1]^2$. The \emph{mesh pattern} $(p,M)$ can be drawn as the permutation graph $G(p) = \{(i,p(i))\}$ in that grid, with the squares of $M$ being shaded. Call $p$ the \emph{pattern} and $M$ the \emph{mesh}. A permutation $w \in S_n$ contains $(p,M)$ if there is an occurrence of $G(p)$ in $G(w)$ (i.e., a $p$-pattern in $w$) in which the regions of $[0,n+1]^2$ corresponding to the mesh in $G(p)$ contain no points of $G(w)$.
\end{definition}

Note that containment of a mesh pattern $(p,\emptyset)$ is equivalent to containment of the classical pattern $p$.

\begin{example}
    Consider the mesh pattern 
    $$\mu = \vcenter
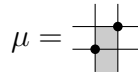
{\hbox{\begin{tikzpicture}[scale=.3]
    \fill[black!20] (1,0) rectangle (2,2);
    \foreach \x in {1,2} {
        \fill (\x,\x) circle (6pt);
        \draw (\x,0) -- (\x,3); 
        \draw (0,\x) -- (3,\x);}
    \end{tikzpicture}}}$$
    whose underlying pattern is $p = 12$. The permutation $1423$ contains four (circled) occurrences of $p$
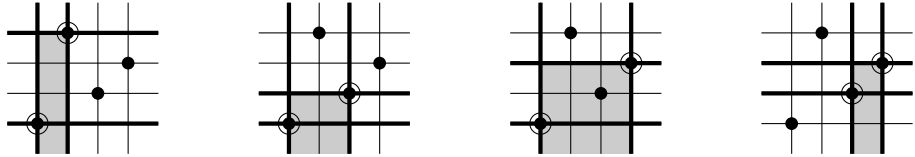
    $$\begin{tikzpicture}[scale=.4]
    \fill[black!20] (1,0) rectangle (2,4);
    \foreach \x in {(1,1), (2,4), (3,2), (4,3)} {
        \fill \x circle (6pt);
    }
    \foreach \x in {(1,1), (2,4)} {
        \draw \x circle (10pt);
    }
    \foreach \x in {1,2} {\draw[ultra thick] (\x,0) -- (\x,5);}
    \foreach \y in {1,4} {\draw[ultra thick] (0,\y) -- (5,\y);}
    \foreach \x in {1,2,3,4} {
        \draw (\x,0) -- (\x,5);
        \draw (0,\x) -- (5,\x);
    }
    \end{tikzpicture}
            \hspace{.5in}
    \begin{tikzpicture}[scale=.4]
    \fill[black!20] (1,0) rectangle (3,2);
    \foreach \x in {(1,1), (2,4), (3,2), (4,3)} {
        \fill \x circle (6pt);
    }
    \foreach \x in {(1,1), (3,2)} {
        \draw \x circle (10pt);
    }
    \foreach \x in {1,3} {\draw[ultra thick] (\x,0) -- (\x,5);}
    \foreach \y in {1,2} {\draw[ultra thick] (0,\y) -- (5,\y);}
    \foreach \x in {1,2,3,4} {
        \draw (\x,0) -- (\x,5);
        \draw (0,\x) -- (5,\x);
    }
    \end{tikzpicture}
            \hspace{.5in}
    \begin{tikzpicture}[scale=.4]
    \fill[black!20] (1,0) rectangle (4,3);
    \foreach \x in {(1,1), (2,4), (3,2), (4,3)} {
        \fill \x circle (6pt);
    }
    \foreach \x in {(1,1), (4,3)} {
        \draw \x circle (10pt);
    }
    \foreach \x in {1,4} {\draw[ultra thick] (\x,0) -- (\x,5);}
    \foreach \y in {1,3} {\draw[ultra thick] (0,\y) -- (5,\y);}
    \foreach \x in {1,2,3,4} {
        \draw (\x,0) -- (\x,5);
        \draw (0,\x) -- (5,\x);
    }
    \end{tikzpicture}
            \hspace{.5in}
    \begin{tikzpicture}[scale=.4]
    \fill[black!20] (3,0) rectangle (4,3);
    \foreach \x in {(1,1), (2,4), (3,2), (4,3)} {
        \fill \x circle (6pt);
    }
    \foreach \x in {(3,2), (4,3)} {
        \draw \x circle (10pt);
    }
    \foreach \x in {3,4} {\draw[ultra thick] (\x,0) -- (\x,5);}
    \foreach \y in {2,3} {\draw[ultra thick] (0,\y) -- (5,\y);}
    \foreach \x in {1,2,3,4} {
        \draw (\x,0) -- (\x,5);
        \draw (0,\x) -- (5,\x);
    }
    \end{tikzpicture}
$$
    but only three occurrences of the mesh pattern $\mu$, because the third picture above containing a point (namely, $(3,w(3)) = (3,2)$) inside the area corresponding to the mesh.
\end{example}

\begin{theorem}\label{thm:barred patterns in P}
A covering relation $v \lessdot w$ exists in the poset $\poset_n$ if and only if 
$v$ and $w$ differ exactly by an occurrence of the mesh pattern
$\vcenter{\hbox{\begin{tikzpicture}[scale=.3]
    \fill[black!30] (1,0) rectangle (2,2);
    \foreach \x in {1,2} {
        \fill (\x,\x) circle (6pt);
        \draw (\x,0) -- (\x,3); 
        \draw (0,\x) -- (3,\x);}
    \end{tikzpicture}}}$
in $v$ becoming an occurrence of 
$\vcenter{\hbox{\begin{tikzpicture}[scale=.3]
    \fill[black!30] (1,0) rectangle (2,2);
    \foreach \x in {1,2} {
        \fill (\x,3-\x) circle (6pt);
        \draw (\x,0) -- (\x,3); 
        \draw (0,\x) -- (3,\x);}
    \end{tikzpicture}}}$
in $w$.
\end{theorem}

\begin{proof}
Suppose that $v \lessdot w$ in $\poset_n$. Then, for some $i$, we have $\inv(v) = (x_1, \ldots, x_n)$ and $\inv(w) = (x_1, \ldots, x_{i-1}, x_i+1, x_{i+1}, \ldots, x_n)$. This implies that $x_i <i-1$, since $\inv(w)$ is a valid inversion sequence. If $i$ 
were to be a left-to-right minimum of $v$ (defined as a value $i$ such that no value $j<i$ appears at a previous position), 
then all $j<i$ would appear to the right of $i$ in $v$, contradicting the fact that $x_i <i-1$. 
Thus the entry $i$ must not be a left-to-right minimum. We can therefore define $j$ as the rightmost entry to the left of $i$ in $v$ such that $j<i$. Then $(j,i)$ is an occurrence of the mesh pattern 
$\vcenter{\hbox{\begin{tikzpicture}[scale=.2]
    \fill[black!30] (1,0.2) rectangle (2,2);
    \foreach \x in {1,2} {
        \fill (\x,\x) circle (7pt);
        \draw (\x,0.2) -- (\x,2.8); 
        \draw (0.2,\x) -- (2.8,\x);}
    \end{tikzpicture}}}$
in $v$. It is then easy to check that the inversion sequence of the permutation obtained by swapping $i$ and $j$ in $v$ is $\inv(w)$. Since the map $\inv$ is a bijection, the permutation obtained by swapping those values in $v$ must be exactly the permutation $w$, proving one direction of the equivalence.

Conversely, assume that $v$ and $w$ differ only by an occurrence of the mesh pattern
$\vcenter{\hbox{\begin{tikzpicture}[scale=.2]
    \fill[black!30] (1,0) rectangle (2,2);
    \foreach \x in {1,2} {
        \fill (\x,\x) circle (6pt);
        \draw (\x,0) -- (\x,3); 
        \draw (0,\x) -- (3,\x);}
    \end{tikzpicture}}}$
in $v$ becoming an occurrence of 
$\vcenter{\hbox{\begin{tikzpicture}[scale=.2]
    \fill[black!30] (1,0) rectangle (2,2);
    \foreach \x in {1,2} {
        \fill (\x,3-\x) circle (6pt);
        \draw (\x,0) -- (\x,3); 
        \draw (0,\x) -- (3,\x);}
    \end{tikzpicture}}}$
in $w$. Denote by $i$ the value of the leftmost entry of this occurrence in $w$. We easily check that $I(w)$ and $I(v)$ differ only at their $i$th coordinate, and that the value at this coordinate is one more in $I(w)$ than in $I(v)$, concluding the proof. 
\end{proof}

We draw the attention of the reader to characterizations of the weak orders and the Bruhat orders in terms of mesh patterns, which are analogous to the statement of Theorem~\ref{thm:barred patterns in P} for the middle order. 
Specifically, it can  be proved that 
\begin{itemize}
    \item any covering relation in the weak order corresponds to an occurrence of 
$\vcenter{\hbox{\begin{tikzpicture}[scale=.2]
    \fill[black!30] (1,0.2) rectangle (2,2.8);
    \foreach \x in {1,2} {
        \fill (\x,\x) circle (7pt);
        \draw (\x,0.2) -- (\x,2.8); 
        \draw (0.2,\x) -- (2.8,\x);}
    \end{tikzpicture}}}$
becoming an occurrence of 
$\vcenter{\hbox{\begin{tikzpicture}[scale=.2]
    \fill[black!30] (1,0.2) rectangle (2,2.8);
    \foreach \x in {1,2} {
        \fill (\x,3-\x) circle (7pt);
        \draw (\x,0.2) -- (\x,2.8); 
        \draw (0.2,\x) -- (2.8,\x);}
    \end{tikzpicture}}}$,
    \item any covering relation in the Bruhat order corresponds to an occurrence of $\vcenter{\hbox{\begin{tikzpicture}[scale=.2]
    \fill[black!30] (1,1) rectangle (2,2);
    \foreach \x in {1,2} {
        \fill (\x,\x) circle (7pt);
        \draw (\x,0.2) -- (\x,2.8); 
        \draw (0.2,\x) -- (2.8,\x);}
    \end{tikzpicture}}}$
becoming an occurrence of 
$\vcenter{\hbox{\begin{tikzpicture}[scale=.2]
    \fill[black!30] (1,1) rectangle (2,2);
    \foreach \x in {1,2} {
        \fill (\x,3-\x) circle (7pt);
        \draw (\x,0.2) -- (\x,2.8); 
        \draw (0.2,\x) -- (2.8,\x);}
    \end{tikzpicture}}}$. 
\end{itemize}

With this perspective, there is an obvious intermediate ordering between the weak and Bruhat orders, and Theorem~\ref{thm:barred patterns in P} shows that it is exactly the poset $\poset_n$. This also explains the name ``middle order'' that we have chosen for $\poset_n$.

Of course, another intermediate ordering could be defined, whose covering relations correspond to an occurrence of
$\vcenter{\hbox{\begin{tikzpicture}[scale=.2]
    \fill[black!30] (1,1) rectangle (2,2.8);
    \foreach \x in {1,2} {
        \fill (\x,\x) circle (7pt);
        \draw (\x,0.2) -- (\x,2.8); 
        \draw (0.2,\x) -- (2.8,\x);}
    \end{tikzpicture}}}$
becoming an occurrence of 
$\vcenter{\hbox{\begin{tikzpicture}[scale=.2]
    \fill[black!30] (1,1) rectangle (2,2.8);
    \foreach \x in {1,2} {
        \fill (\x,3-\x) circle (7pt);
        \draw (\x,0.2) -- (\x,2.8); 
        \draw (0.2,\x) -- (2.8,\x);}
    \end{tikzpicture}}}$.
This yields a poset isomorphic to the middle order. 
This isomorphic version can also be described from the same order on inversion sequences, however associating a permutation to an inversion sequence in a different way. 
More precisely, it follows from building the inversion sequence $(z_1, \dots, z_n)$ associated with a permutation $w$ of size $n$ by computing the number $z_i$ of inversions in which $n-i$ is the inversion bottom.

Theorem~\ref{thm:barred patterns in P} and the above remark also imply another relationship between the middle order and the Bruhat and weak orders. 

\begin{corollary}\label{cor:P on 231 avoiders}
When restricted to $213$-avoiding permutations, the middle order is identical to the Bruhat order.
\end{corollary}

\begin{corollary}\label{cor:P on 132 avoiders}
When restricted to $132$-avoiding permutations, the middle order is identical to the weak order.
\end{corollary}

\section{Connections with previous works}\label{sec:ipf}

The middle order on permutations has occasionally surfaced in the literature, although not under this name. The present section aims at giving a (probably not exhaustive) survey of such appearances.

To the best of our knowledge, the first appearance of the middle order, actually in a much more general context, occurred in \cite{A}, where Armstrong defines a new partial order associated to each element of an arbitrary Coxeter group (the middle order corresponds to Armstrong's sorting order based on the long element in the symmetric group). 
A few years later, an isomorphic version of the middle order was defined in \cite{D}, where some basic combinatorial properties are proved. Subsequently, a far reaching generalization was introduced in \cite{CG}, in terms of the notion of a \emph{$\delta$-cliff}, and some general algebraic properties are shown for the resulting $\delta$-cliff posets. 

More recently, in \cite{CDMY}, the authors consider ``interval parking functions,'' which are defined like classical parking functions with the constraint that each car has an interval of preferred spots (instead of having a single preferred spot to park in). Thus, in order to describe an interval parking preference, one can do so by giving two classical parking preferences $a=(a(1),\ldots ,a(n))$ and $b=(b(1),\ldots ,b(n))$ such that car $i$ is willing to park only in the interval $[a(i),b(i)]$. The pair $(a,b)$ is an \emph{interval parking function} when all cars can successfully park in their preferred intervals. It can be proved that interval parking functions can be encoded by certain pairs of permutations (of the same size). The set of such pairs actually defines a partial order, called \emph{pseudoreachability order} by the authors, which turns out to be isomorphic to the middle order. Our Theorem~\ref{thm:barred patterns in P} completely answers the question (left open in \cite{CDMY}) of characterizing covering relations in the pseudoreachability order in terms of patterns.

The middle order has also surfaced on several occasions in a probabilistic setting. 
One such recent appearance is in \cite{C}: the author introduces \emph{Mallows processes}, and proves that 
the transition graph of a regular Mallows process is a so-called \emph{expanded hypercube}. 
In our language, the definition of expanded hypercubes from~\cite{C} says that they are the unoriented graphs underlying the Hasse diagrams of the middle orders. 
Staying in the probability literature, \cite{C} is actually related to
the earlier paper \cite{AP} (which studies the probability that a permutation chosen uniformly at random with respect to size and number of inversions is indecomposable). Indeed, \cite{AP} considers Markov processes on permutations that are similar to the processes of~\cite{C}, however with a different purpose. 
Relevant to us is that, in the processes of~\cite{AP}, permutations evolve exactly along the covering relations of our middle order. These are briefly defined in~\cite{AP}, although without a specific name. 
A third appearance in a probabilistic context is in \cite{JLL}, under the name of \emph{grid order}. There, the authors are interested in correlation properties of random permutations analogous to the Harris-Kleitman inequality for up-sets. In this paper, the fact that the middle order ``sits'' between the weak order and the Bruhat order is key, as is the distributive lattice property of the middle order. 

Finally, restricting the middle order to notable sets of permutations gives rise to interesting structures. We study its restriction to involutions in \cref{sec:involutions}, and the forthcoming article~\cite{Nick} explores its restriction to alternating permutations and to 132-avoiding permutations, both of which are related to the study of ``mixed dimer covers'' (which are generalizations of perfect matchings) on some graphs. 

\section{Characterizing and counting intervals}\label{sec:counting}

An important aspect of investigating the structure of a (combinatorially defined) poset is to understand its intervals. This approach has been undertaken in many cases, for instance for Tamari lattices \cite{CCP}, the Bruhat order \cite{T}, the weak order \cite{EHKM}, the consecutive pattern poset \cite{EM}, the matching pattern poset \cite{CF}, and the Dyck pattern poset \cite{BCFS}, to cite just a few. In each setting, these efforts have given insight into the overall architecture and complexity of the poset.

The proof of the next proposition is elementary, 
but we did not find it in the literature and so we include it here. We write $\Int(\mathcal{Q})$ for the set of intervals of a poset $\mathcal{Q}$.

\begin{proposition}\label{prop:intervals in product posets}
    Let $\mathcal{Q}_1,\ldots \mathcal{Q}_h$ be a sequence of finite posets. Then $\Int(\mathcal{Q}_1 \times \cdots \times \mathcal{Q}_h )= \Int(\mathcal{Q}_1)\times \cdots \times \Int(\mathcal{Q}_h)$, and hence
    $$
    |\Int(\mathcal{Q}_1 \times \cdots \times \mathcal{Q}_h) |=\prod_{i=1}^{h}|\Int(\mathcal{Q}_i )|.
    $$
\end{proposition}

\begin{proof}
    Each interval in $\mathcal{Q}_1\times \cdots \times \mathcal{Q}_h$ has the form $[(x_1,\ldots,x_h), (y_1,\ldots, y_h)]$, where $x_i \le y_i$ in $\mathcal{Q}_i$ for all $ 1 \leq i \leq h$. Then 
    \begin{align*}
        [(x_1,\ldots,x_h), \ (y_1,\ldots, y_h)] &= \{(z_1,\ldots, z_h)\, | \, (x_1,\ldots, x_h) \le (z_1, \ldots, z_h) \le (y_1, \ldots, y_h)\}\\
        &= \{(z_1,\ldots,z_h)\, | \, x_i \le z_i \le y_i \text{ for all }i\}\\
        &= [x_1,y_1] \times \cdots \times [x_h,y_h]. \qedhere
    \end{align*}
\end{proof}

This allows us to compute the number of intervals in $\poset_n$.

\begin{lemma}\label{lem:nb_intervals_chain}
    The number of intervals of the $n$-element chain $[0,n-1]$ is $\binom{n+1}{2}$.  
\end{lemma}

\begin{proof}
    Intervals of the $n$-element chain poset $[0,n-1]$ are size-2 multisets of $\{0,1,\ldots, n-1\}$. The number of these is
    \begin{align*} n + \binom{n}{2} &= \binom{n+1}{2}. \qedhere \end{align*}
\end{proof}

From here, the enumeration of intervals in $\poset_n$ follows immediately from  Proposition~\ref{prop:poset bijection},  Proposition~\ref{prop:intervals in product posets} and Lemma~\ref{lem:nb_intervals_chain}.

\begin{corollary}\label{cor:counting intervals in P}
    The number of intervals in the poset $\poset_n$ is 
    $$\prod_{i=0}^{n-1}|\Int([0,i])| = \prod_{i=0}^{n-1} \binom{i+2}{2} =\frac{n!(n+1)!}{2^n}.$$ 
\end{corollary}

As intervals of a graded poset, the intervals of $\poset_n$ are also graded posets. We next show that their enumeration can be refined by rank.

\begin{theorem}\label{thm:intervals by rank}
    Denote by $f(n,k)$ the number of intervals in $\poset_n$ having rank $k$, with $n\geq 1$ and $k\geq 0$. Then $f(1,0)=1$ and for $n\geq 2$ and $k \in [0, \binom{n}{2}]$, 
    $$
    f(n,k)=\sum_{h=0}^{n-1}(n-h) \cdot f(n-1,k-h),
    $$
(with the convention that $f(n,j) =0 $ when $j<0$). 
\end{theorem}

\begin{proof}
Again, we shall rely on the isomorphism identified in Proposition~\ref{prop:poset bijection} and prove the announced result for the poset $[0,0] \times [0,1] \times [0,2] \times \cdots \times [0,n-1]$. 
Our proof is  by induction on $n$, with small cases being easy to check. 
By Proposition~\ref{prop:intervals in product posets}, an interval $I$ of $[0,0] \times [0,1] \times [0,2] \times \cdots \times [0,n-1]$ is uniquely characterized as a pair $(J,K)$ consisting of an interval $J$ of $[0,0] \times [0,1] \times [0,2] \times \cdots \times [0,n-2]$ and an interval $K$ of $[0,n-1]$. 
If we denote by $k \in [0,\binom{n}{2}]$ the rank of $I$, and $h \in [0,n-1]$ the rank of $K$, then the rank of $J$ has to be $k-h$. 
For any fixed value of $h$, there are $n-h$ intervals of rank $h$ in $[0,n-1]$ (since those are of the form $[i,i+h] \subseteq [0,n-1]$). 
Similarly, when $h$ is fixed, and assuming $k-h\geq 0$, there are $f(n-1,k-h)$ intervals of rank $k-h$ in $[0,0] \times [0,1] \times [0,2] \times \cdots \times [0,n-2]$ by the inductive hypothesis, yielding the desired result.
\end{proof}

The first few values of $f(n,k)$ are recorded in Table~\ref{table:intervals by rank} below. 

\begin{table}[htbp]
\begin{center}
\begin{tabular}{c|ccccccccccc}
		\backslashbox{$n$}{$k$} & 0 & 1 & 2 & 3 & 4 & 5 & 6 & 7 & 8 & 9 & 10 \\
		\hline
		1 & 1   &     &     &      &    &     &     &    &    &   &    \\
		2 & 2   & 1   &     &      &    &     &     &    &    &   &    \\
		3 & 6   & 7   & 4   & 1    &    &     &     &    &    &   &    \\
		4 & 24  & 46  & 49  & 36  & 18  & 6   & 1   &    &    &   &    \\
		5 & 120 & 326 & 501 & 562 & 497 & 354 & 204 & 94 & 33 & 8 & 1  \\
		
	\end{tabular}
	\end{center}
 \bigskip
	\caption{The number of intervals in $\poset_n$ having rank $k$.}\label{table:intervals by rank}
\end{table}

For $k=0$, we are just counting elements of $\poset_n$, thus $f(n,0)=n!$.

For $k=1$ we are counting covering relations in $\poset_n$. The corresponding sequence $f(n,1)=n!(n-H_n)$, where $H_n=\sum_{i=1}^{n}\frac{1}{i}$ is the $n$th harmonic number, is \cite[A067318]{oeis}, as it can be shown by specializing the recurrence given in Theorem~\ref{thm:intervals by rank} for $k=1$. This turns out also to equal the sum of reflection lengths of all elements of $S_n$, 
but we have not found a bijective explanation of this fact yet. 

The full table of data appears in \cite[A139769]{oeis}. That entry refers to ``term-like derivatives,'' and the coincidence of this data follows from our recursion and the fact that the row polynomials are also obtained recursively. 
More precisely, the coefficients $p(n,k)$ of the polynomials $p_n (x)=\sum_{k=0}^{\binom{n}{2}}p(n,k)x^k$ defined by $p_n (x) =\prod_{i=1}^{n}\sum_{j=1}^{i}D(x^j )$ (where $D$ is the usual derivative operator), for $n\geq 1$, appear as the rows of table A139769 (when rows are indexed starting from 1 and columns are indexed starting from 0). It is clear that $p_n (x)=p_{n-1}(x)\cdot \sum_{j=1}^{n}D(x^j )$. This polynomial equality turns out to be equivalent to the recurrence proved in Theorem \ref{thm:intervals by rank} (upon a reversal of the rows of Table \ref{table:intervals by rank}). In fact, we have that
\[
\sum_{k=0}^{\binom{n}{2}}p(n,k)x^k =\left( \sum_{k=0}^{\binom{n-1}{2}}p(n-1,k)x^k \right) \cdot \left( \sum_{k=0}^{n-1}(k+1)x^k \right) ,
\]
from which we get
\[
p(n,k)=\sum_{h=0}^{k}(h+1)\cdot p(n-1,k-h).
\]

Now, if we define $q_n (x)=x^{\binom{n}{2}}\cdot p_n (\frac{1}{x})$, and write its coefficients as $q_n (x)=\sum\limits_{k=0}^{\binom{n}{2}}q(n,k)x^k $, the above equation yields 

\begin{align*}
q(n,k)&=p\left( n,\binom{n}{2}-k\right) =\sum_{h=0}^{\binom{n}{2}-k}(h+1)\cdot p\left( n-1,\binom{n}{2}-k-h \right) \\
&=\sum_{h=0}^{\binom{n}{2}-k}(h+1)\cdot q\left( n-1,\binom{n-1}{2}-\left( \binom{n}{2}-k-h\right) \right) \\
&=\sum_{h=0}^{\binom{n}{2}-k}(h+1)\cdot q(n-1,(k+h)-(n-1)) =\sum_{h=k-\binom{n-1}{2}}^{n-1}(n-h)\cdot q(n-1,k-h).
\end{align*}
Observing that $q(n-1,k-h)=0$ when $0 \leq h< k-\binom{n-1}{2}$, this is exactly our recurrence for $f(n,k)$ in Theorem \ref{thm:intervals by rank}.

\bigskip

In the rest of this section we restrict our attention to a particularly nice family of intervals. Specifically, we provide a complete characterization and enumeration of boolean intervals in $\poset_n$. 
Boolean intervals are among the nicest substructures one can find in a poset. For instance, and this is further discussed at the end of the present section,  the computation of the M\"obius function of boolean intervals is easy. This gives important information on the topological properties of a poset. Let us mention that boolean intervals have been studied in the Bruhat order in \cite{T, T-intervalstructures}, and in the weak order in \cite{EHKM, T-intervalstructures}. 

\begin{definition}
    Let $J$ be an interval in a poset $P$. If $J$ is isomorphic to a boolean algebra, then $J$ is a \emph{boolean} interval.
\end{definition}

By Proposition~\ref{prop:poset bijection}, an interval in $\poset_n$ corresponds to a pair of inversion sequences 
$$(x_1, \ldots, x_n) \le (y_1, \ldots, y_n)$$
with $x_i \le y_i$ for all $i$. This leads to an immediate characterization of the boolean intervals in $\poset_n$.

\begin{lemma}\label{lem:boolean interval characterization}
An interval $[v,w] \subseteq \poset_n$ is 
boolean if and only if, for $\inv(v) = (x_1,\ldots, x_n)$ and $\inv(w) = (y_1,\ldots,  y_n)$, we have $y_i \in \{x_i, x_i+1\}$ for all $i$. 

Moreover, the boolean interval $[v,w]$ has rank $k$ if and only if there are exactly $k$ indices $i$ such that $y_i = x_i+1$.
As a consequence, the rank of a boolean interval in $\poset_n$ is at most $n-1$. 
\end{lemma}

\begin{proof}
    The interval defined by the two elements $v$ and $w$ is isomorphic to the product of intervals
    $$[x_1,y_1]\times \cdots \times [x_n,y_n],$$
    and a product of chains is boolean if and only if each of those chains has rank $0$ or $1$. 
    The rank property follows since rank $k$ is obtained when choosing exactly $k$ of these chains to be of rank $1$.
    The final statement follows from the requirement that $x_1 = y_1 = 0$.
\end{proof}

In analogy to Corollary~\ref{cor:counting intervals in P}, which counted all intervals in $\poset_n$, we can also enumerate the (nicely behaved) boolean intervals of $\poset_n$. The sequence described in Corollary~\ref{cor:boolean interval count} appears as \cite[A001147]{oeis}. 

\begin{corollary}\label{cor:boolean interval count}
    The number of boolean intervals in $\poset_n$ is $(2n-1)!!$.
\end{corollary}

\begin{proof}
With the notation of Lemma~\ref{lem:boolean interval characterization}, for each $i$, the pair $(x_i,y_i)$ has $i$ possibilities if $y_i = x_i$, and $i-1$ possibilities if $y_i = x_i+1$. Thus there are $2i-1$ total choices, for each value of $i$.
\end{proof}

For $n > k \ge 0$, we denote by $b(n,k)$ the number of rank-$k$ boolean intervals in $\poset_n$. 
These values are recorded for $n \le 5$ in Table~\ref{table:boolean intervals by rank}. 
Those values match \cite[A059364]{oeis}. That sequence is defined by the closed formula $d(n,k)=\sum_{i=0}^{n}\binom{i}{k}c(n,n-i)$, where the $c(n,j)$ are the signless Stirling numbers of the first kind. In other words, for all $n\geq j \geq 0$,
\begin{equation}\label{equation:signless stirling}
    c(n,j) \text{ is the number of permutations of size $n$ having $j$ cycles. }
\end{equation}
The next theorem proves that $d(n,k)$ is indeed equal to our $b(n,k)$.

\begin{table}[htbp]
\begin{center}
	\begin{tabular}{c|ccccc}
		\backslashbox{$n$}{$k$} & 0 & 1 & 2 & 3 & 4 \\
		\hline
		1 & 1   &     &     &     &    \\
		2 & 2   & 1   &     &     &    \\
		3 & 6   & 7   & 2   &     &    \\
		4 & 24  & 46  & 29  & 6   &    \\
		5 & 120 & 326 & 329 & 146 & 24 \\
		
	\end{tabular}
	\end{center}
 \bigskip
	\caption{The number $b(n,k)$ of boolean intervals in $\poset_n$ having rank $k$, for small $n$ and $k$.}\label{table:boolean intervals by rank}
\end{table}

\begin{theorem}\label{thm:formula for b(n,k)}
Fix $n > k \ge 0$. The number of boolean intervals in $\poset_n$ having rank $k$ is 
$$b(n,k)=\sum_{i=0}^{n}\binom{i}{k}c(n,n-i).$$
\end{theorem}

\begin{proof}
We use a slight modification of Foata's fundamental bijection, defined as follows. Given a permutation $w$, write $w$ in cycle notation such that, in each cycle, the smallest element is listed last and the cycles are sorted in increasing order of their minima. Let $s$ be the word obtained by deleting the parentheses from this expression. The image of $w$ under the bijection is the permutation whose one-line notation is $s$.  

A consequence of this bijection is that $c(n,j)$ also counts the permutations of size $n$ having $j$ right-to-left minima (defined as values $i$ such that no value $j<i$ appears at a later position). 
Also, by definition of $I$, for a permutation $w$ of size $n$, with inversion sequence $I(w) = (x_1, \dots, x_n)$, we have $x_i=0$ if and only if the value $i$
is a right-to-left minimum of $w$. 
It follows that for $i \in [0,n-1]$, 
the value $c(n,n-i)$ counts inversion sequences of size $n$ having $i$ non-zero entries. 

We saw in \cref{lem:boolean interval characterization} that a boolean interval of rank $k$ in  $\poset_n$ corresponds to a pair of inversion sequences $x=(x_1,\ldots, x_n)$ and $y=(y_1,\ldots,  y_n)$ such that $y_i=x_i+1$ for exactly $k$ indices $i$, and $y_i=x_i$ for all other indices $i$. 
Such a boolean interval is therefore described by an inversion sequence $(y_1,\ldots,  y_n)$ in which $k$ non-zero entries are marked (the sequence $x$ being then obtained from $y$ by subtracting $1$ from each marked entry). Since there are $\binom{i}{k}$ ways of choosing these non-zero entries when $y$ has $i$ non-zero entries, we obtain $b(n,k)=\sum_{i=0}^{n}\binom{i}{k}c(n,n-i)$, as desired. 
\end{proof}

We note that the OEIS indicates that the values $b(n,k)$ also obey the recursive relation 
$$b(n,k)=n \cdot b(n-1,k)+(n-1) \cdot b(n-1,k-1).$$
This can be easily proved from the combinatorial characterization of boolean intervals given in Lemma~\ref{lem:boolean interval characterization}.

As mentioned above, boolean intervals in a poset interact neatly with a particularly important function in the study of posets, namely the M\"obius function.

\begin{definition}\label{defn:mobius}
    The \emph{M\"obius function} $\mu$ on a poset $P$ is defined recursively by
    $$\mu(s,u) = \begin{cases}
        0 & \text{if $s \not\le u$},\\
        1 & \text{if }s = u, \text{ and}\\
        -\!\!\sum\limits_{s \le t < u} \mu(s,t) & \text{for all $s < u$}.
    \end{cases}$$
\end{definition}

Since $\poset_n$ is a finite distributive lattice, it is well-known (see, for example, \cite{ec1}) that, for any $v,w\in \poset_n$, the M\"obius function $\mu (v,w)$ is equal to 0 if the interval $[v,w]$ is not boolean, and otherwise $\mu (v,w)=(-1)^t$, where $t$ is the rank of $[v,w]$. 
So, for instance, if we look at principal order ideals in $\poset_n$, writing $\id$ for the identity permutation and taking $w\neq \id$, we have that $\mu (\id ,w)=0$ if and only if $I(w)$ contains at least one entry that is at least $2$ (equivalent to $w$ containing a $312$- or $321$-pattern). Otherwise, $\mu (\id ,w) =(-1)^t$, where $t$ is the number of 1's in $I(w)$ (which, in this case, is the total number of inversions of $w$).

\section{The Euler characteristic}\label{euler_char}

We have seen in Proposition~\ref{prop:poset bijection} that the middle order posets $\poset_n$ are isomorphic to direct products of chains, and hence are finite distributive lattices. 
In fact, they are very special among finite distributive lattices. 
Indeed, every finite distributive lattice is isomorphic to a sublattice of a product of chains. This may be seen as a consequence of Birkhoff's representation theorem for finite distributive lattices \cite{B}, asserting that every finite distributive lattice $L$ is isomorphic to a sublattice of the boolean algebra $\wp (L)$ (i.e., of the product of chains of size 2), even if this is certainly not the \emph{minimal} product of chains into which $L$ can be embedded in general.

It is then reasonable to study middle order posets from the perspective of distributive lattices. With this in mind, in the present section we will explore a particular combinatorial aspect of those objects: the Euler characteristic. We start by giving some necessary preliminary definitions and results.

In a finite distributive lattice $L$, we will write $\0$ for the minimum element of $L$. An element $a\neq \0$ of $L$ is said to be \emph{join-irreducible} when the only way to write $a$ as the join of two elements is if one of those elements is $a$ itself. For finite lattices (like each $\poset_n$), this is equivalent to saying that $a$ covers exactly one element. 

It is easy to characterize join-irreducible elements of $\poset_n$.

\begin{lemma}
\label{lem:join_irred}
A permutation $w$ of size $n$ is join-irreducible in $\poset_n$ if and only if it is written in one-line notation as 
\[
w = 1\ 2\ \cdots \ i \ j\ (i+1)\ (i+2)\ \cdots \ (j-1)\ (j+1)\ \cdots\ n,
\]
for some $i \in[0,n-2]$ and $j \in [i+2,n]$. 
\end{lemma}

\begin{proof}
Because join-irreducible elements are those covering exactly one element, a permutation $w$ is join-irreducible exactly when its inversion sequence $\inv(w)$ has exactly one nonzero entry (the only element covered by $w$ then being obtained by subtracting $1$ from this entry). 
It is therefore enough to describe, for each $j \leq n$ and  $k \in [1,j-1]$, the permutation whose inversion sequence has an entry of value $k$ in position $j$, and all other entries equal to $0$. 
This permutation is easily seen to be the one claimed in the lemma, for $i=j-1-k$. 
Indeed, the only inversion top in this permutation is $j$, and $j$ is the top of $k=j-1-i$ inversions. 
\end{proof}

\begin{remark}
We can observe that the cycle representation of the permutation $w$ appearing in Lemma~\ref{lem:join_irred} consists of a single cycle, of size $j-i$, which is precisely $(j\ (j-1) \ \cdots \ (i+2)\ (i+1))$. Therefore, Lemma~\ref{lem:join_irred} can equivalently be stated as follows: a permutation is join-irreducible in $\poset_n$ if and only if its cycle representation consists of a single cycle of size at least 2, which can be written so that its elements appear in decreasing order, and those elements are consecutive in value.
\end{remark}

\begin{definition}\label{defn:valuation}
A \emph{valuation} of a finite distributive lattice $L$ (with values in a commutative ring $\mathbf{A}$ with identity) is a map $\val:L\rightarrow \mathbf{A}$ such that $\val(\0)=0$ and, for any $x,y\in L$, 
$$\val(x)+\val(y)=\val(x\wedge y)+\val(x\vee y).$$ 
\end{definition}

It is known that a valuation is uniquely determined by its values on the join-irreducibles, and such values can be arbitrarily assigned \cite{R}.

The \emph{Euler characteristic} of $L$ is defined as the unique valuation $\chi$ such that $\chi (a)=1$ for every join-irreducible $a \in L$.

Valuations for distributive lattices are, in a sense, analogs of measures for boolean algebras, and the Euler characteristic is an example of a valuation that has no counterpart in measure theory. From a combinatorial point of view, the Euler characteristic of a distributive lattice $L$ is connected with the M\"obius function of the poset $J(L)$ of the join-irreducibles of $L$. From a more geometric point of view, it has an important homological interpretation: it is the classical Euler characteristic in the order homology of the poset $J(L)$. 
For information on the order homology of a poset and related topics we refer the reader to \cite{R,ec1}.

\begin{theorem}\label{thm:euler characteristic}
    The Euler characteristic $\chi (w)$ of a permutation $w\in \poset_n$ is equal to the number of right-to-left non-minima of $w$. Moreover, 
    the number of permutations of size $n$ having Euler characteristic $k$ 
    is $c(n,n-k)$, as defined in~\eqref{equation:signless stirling}.  
\end{theorem}

\begin{proof} 
Recall from the proof of Theorem~\ref{thm:formula for b(n,k)} that the right-to-left minima of a permutation correspond to zeroes in its inversion sequence. Define the map $\chi$ by setting $\chi (v)$ equal to the number of right-to-left non-minima of $v$, equivalently, the number of non-zero entries in its inversion sequence. 
Clearly $\chi (\id )=0$. 
It is also easy to see that $\chi(a) = 1$ for all join-irreducible elements $a \in \poset_n$ (since $\inv (a)$ has exactly one nonzero entry for such $a$).
Fix $v,w \in \poset_n$, 
with $\inv(v) = (x_1,\ldots, x_n)$ and $\inv(w) = (y_1,\ldots, y_n)$. By the description of the meet and join operations in $\poset_n$ (see \eqref{eq:meet_and_join1} and \eqref{eq:meet_and_join2}), we see that 
    \begin{align*}
        \chi(v\vee w) &= n - \#\{i\, |\, x_i = y_i = 0\}, \text{ and}\\
        \chi(v\wedge w) &= n - \left(\#\{i\, |\, x_i = 0\} + \#\{i\, |\, y_i = 0\} - \#\{i\, |\, x_i = y_i = 0\}\right)\\
        &= \chi(v) + \chi(w) - n + \#\{i\, |\, x_i = y_i = 0\}.
    \end{align*}
    Thus $\chi$ is a valuation, and therefore is the Euler characteristic. 
    Moreover, recalling that the signless Stirling number of the first kind $c(n,j)$ counts permutations of size $n$ having $j$ cycles, by Foata's bijection (or rather, its modified version used in the proof of Theorem~\ref{thm:formula for b(n,k)}) we have that $c(n,n-k)$ is also the number of permutations of size $n$ having $k$ right-to-left non-minima, 
    as desired.  
\end{proof}

\section{The poset of involutions and its M\"obius function}\label{sec:involutions}

Involutions form a natural subset of permutations, and subposets of involutions have been investigated in the case of both the Bruhat order \cite{I,CH} and the weak order \cite{CJW}. For instance, the posets of involutions under the Bruhat order have many desirable properties: they are graded, lexicographically shellable (hence Cohen-Macaulay), and Eulerian. 

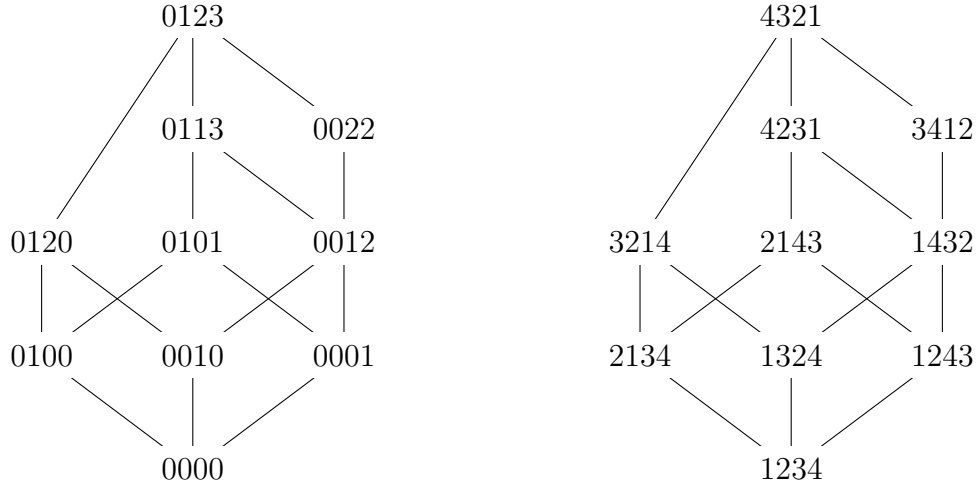
\begin{figure}[htbp]
\begin{center}
\begin{tikzpicture}
    \draw (0,0) coordinate (0000);
    \draw (2,1.5) coordinate (0001);
    \draw (0,1.5) coordinate (0010);
    \draw (-2,1.5) coordinate (0100);
    \draw (2,3) coordinate (0012);
    \draw (0,3) coordinate (0101);
    \draw (-2,3) coordinate (0120);
    \draw (2,4.5) coordinate (0022);
    \draw (0,4.5) coordinate (0113);
    \draw (0,6) coordinate (0123);
    \foreach \x in {0100,0010,0001} {\draw (0000) -- (\x);}
    \foreach \y in {0120, 0101} {\draw (0100) -- (\y);}
    \foreach \y in {0120, 0012} {\draw (0010) -- (\y);}
    \foreach \y in {0101, 0012} {\draw (0001) -- (\y);}
    \foreach \x in {0101,0012,0123} {\draw (0113) -- (\x);}
    \draw (0022) -- (0012);
    \foreach \x in {0120,0022} {\draw (0123) -- (\x);}
    \foreach \x in {0000,0100,0010,0001,0120,0101,0012,0113,0022,0123} {\fill[white] (\x)++(-.45,-.3) rectangle ++(.9,.6); \draw (\x) node {$\x$};}
\end{tikzpicture}
    \hspace{1in}
\begin{tikzpicture}
    \draw (0,0) coordinate (1234);
    \draw (2,1.5) coordinate (1243);
    \draw (0,1.5) coordinate (1324);
    \draw (-2,1.5) coordinate (2134);
    \draw (2,3) coordinate (1432);
    \draw (0,3) coordinate (2143);
    \draw (-2,3) coordinate (3214);
    \draw (2,4.5) coordinate (3412);
    \draw (0,4.5) coordinate (4231);
    \draw (0,6) coordinate (4321);
    \foreach \x in {2134,1324,1243} {\draw (1234) -- (\x);}
    \foreach \y in {3214, 2143} {\draw (2134) -- (\y);}
    \foreach \y in {3214, 1432} {\draw (1324) -- (\y);}
    \foreach \y in {2143, 1432} {\draw (1243) -- (\y);}
    \foreach \x in {2143,1432,4321} {\draw (4231) -- (\x);}
    \draw (3412) -- (1432);
    \foreach \x in {3214,3412} {\draw (4321) -- (\x);}
    \foreach \x in {1234,2134,1324,1243,3214,2143,1432,4231,3412,4321} {\fill[white] (\x)++(-.45,-.3) rectangle ++(.9,.6); \draw (\x) node {$\x$};}
\end{tikzpicture}

\end{center}
\caption{The poset $\invol_4$ of involutions of size 4 (on the right) and its isomorphic poset defined by inversion sequences (on the left).}\label{I_4}
\end{figure}

If we restrict the middle order to involutions, we obtain certain subposets $\invol_n$ of our lattices $\poset_n$, which appear to be much more difficult to investigate than the lattices in their entirety. In particular, the $\invol_n$ are not lattices, and they are not even graded, as is evident from Figure~\ref{I_4}. Moreover, the poset $\invol_n$ is not an \emph{interval-closed} subposet of $\poset_n$, meaning that there exist $v,w\in \invol_n$ and $z\in \poset_n \setminus \invol_n$ such that $v\leq z\leq w$. This makes it not obvious at all, for instance, how to compute the M\"obius function of $\invol_n$. Rather surprisingly, it turns out that the M\"obius function of the principal order ideals of $\invol_n$ has a neat and purely combinatorial expression. Presenting this property is the goal of the present section. 

As introduced above, we denote by 
$\invol_n$
the subposet of $\poset_n$ consisting of the involutions, and we set 
$i(n):=|\invol_n|$ the number of involutions of size $n$. 
The first tool that we need is a characterization of involutions in terms of their inversion sequences. To do this, we recall the following well-known recursion satisfied by $i(n)$:
$$
i(n) = \begin{cases}
i(n-1)+(n-1) \cdot i(n-2) & \textnormal{ for $n\geq 2$, and}\\
1 & \text{ for } n \in \{0,1\}. 
\end{cases}
$$

This is easily proved by observing that, for a given involution of size $n$, either $n$ is a fixed point (which gives the term $i(n-1)$) or $n$ belongs to a 2-cycle $(k\, n)$, for some $k\in [1,n-1]$ (which gives the term $(n-1) \cdot i(n-2)$).

This same decomposition of involutions can be used to characterize their inversion sequences as follows.

\begin{proposition}\label{rec_char_inv}
Let $\inv(w)=(x_1 ,\ldots ,x_n)$ be the inversion sequence of a permutation $w$ of size $n$. Then $w$ is an involution if and only if
\begin{itemize}
\item[(i)] $x_n =0$ and $(x_1 ,\ldots x_{n-1})$ is the inversion sequence of an involution of size $n-1$, or
\item[(ii)] $x_n =k>0$, $x_{n-k}=0$ and $(x_1 ,\ldots ,x_{n-k-1},x_{n-k+1}-1,\ldots ,x_{n-1}-1)$ is the inversion sequence of an involution of size $n-2$. 
\end{itemize}
\end{proposition}

\begin{proof}
    First suppose that $w(n) = n$, which is equivalent to assuming that $x_n = 0$. Let $u \in S_{n-1}$ be the permutation defined by restricting $w$ to $[1,n-1]$. 
    This has inversion sequence $\inv(u) = (x_1,\ldots,x_{n-1})$. Then $w \in \invol_n$ if and only if $u \in \invol_{n-1}$; i.e., if and only if 
    $\inv(u)$ is the inversion sequence of an involution.

    Now suppose that $w(n) < n$, which is equivalent to assuming that $x_n = k$ for some $k > 0$; equivalently, $w^{-1}(n) = n-k$. Let $v \in S_{n-2}$ be the permutation obtained by deleting both $n$ and $n-k$ from the one-line notation of $w$ and then ``flattening,'' that is, replacing $n-j$ by $n-j-1$ for any $j\in [1,k-1]$. Then $w \in \invol_n$ if and only if $w(n) = n-k$ and $v \in \invol_{n-2}$. 
    We claim that $w(n) = n-k$ if and only if 
    \begin{itemize}
        \item $x_{n-k} = 0$ and 
        \item all $m > n-k$ appear to the left of $n-k$ in the one-line notation for $w$.
    \end{itemize}
    Indeed, $w(n) = n-k$ means that no entry of the permutation appears to the right of the value $n-k$, while the first condition (resp., the second condition) above indicates that no entry smaller than $n-k$ (resp., larger than $n-k$) appears to the right of the value $n-k$. 
    The second condition can be rewritten as $x_m > 0$ for all $m > n-k$. 
    Thus $w$ is an involution if and only if $x_{n-k} = 0$ and $v$ is an involution with inversion sequence $\inv(v) = (x_1,\ldots, x_{n-k-1}, x_{n-k+1} - 1, \ldots, x_{n-1} - 1)$.
\end{proof}

We need to introduce some terminology for describing sequences. To start, the \emph{slowly increasing sequence of size $n$} is the sequence
$$(0, 1, 2, \ldots, n-1).$$

\begin{definition}\label{defn:ascents and relevance}
    Let $x=(x_1 ,\ldots ,x_n)$ be an inversion sequence of size $n$. If $x_i < x_{i+1}$, then $i$ is an \emph{ascent}. We say that $i$ is a \emph{small ascent} when $x_{i+1}=x_i +1$, and an ascent is \emph{large} when it is not small. 
    We say that $x$ is \emph{slow-climbing} when $x$ does not contain large ascents. An involution $w$ is \emph{slow-climbing} (resp., \emph{not slow-climbing}) when its inversion sequence $I(w)$ is slow-climbing (resp., not slow-climbing). 
\end{definition}

The slow-climbing quality (or not) of an inversion sequence will be critical to our calculation of the M\"obius function in these posets of involutions.

\begin{lemma} \label{lem: relevant = slowly increasing}
The inversion sequence of an involution is slow-climbing if and only if it is the concatenation of slowly increasing sequences.  
\end{lemma}

\begin{proof}
	It is easy to check that any concatenation of slowly increasing sequences is the inversion sequence of an involution (specifically, of the direct sum of reverse identity permutations), and clearly such a sequence is slow-climbing.
	
	Conversely, suppose that $x$ is a slow-climbing inversion sequence of an involution. We proceed by induction on the size $n$ of $x$. We can use Proposition~\ref{rec_char_inv} and consider two distinct cases. If $x_n=0$, then $(x_1,\ldots ,x_{n-1})$ is the inversion sequence of an involution as well, and by induction it is the concatenation of slowly increasing sequences. Therefore the same is clearly true for $x=(x_1 ,\ldots, x_{n-1},0)$. Now suppose that $x_n =k>0$. Then, by	Proposition~\ref{rec_char_inv}, $x_{n-k}=0$ and $(x_1 ,\ldots ,x_{n-k-1},x_{n-k+1}-1,\ldots, x_{n-1}-1)$ is the inversion sequence of an involution of size $n-2<n$. 
 Moreover, since $x$ is slow-climbing and satisfies $x_{n-k}=0$ and $x_n=k$, then necessarily $x_{n-k+i}=i$ for all $0\leq i\leq k$; i.e., $(x_{n-k},\ldots, x_n)$ is a slowly increasing sequence. 
 Therefore $(x_1 ,\ldots ,x_{n-k-1},x_{n-k+1}-1,\ldots, x_{n-1}-1)$ is slow-climbing, and by induction it must be that $(x_1 ,\ldots ,x_{n-k-1},x_{n-k+1}-1,\ldots, x_{n-1}-1)$ is the concatenation of slowly increasing sequences, where in particular $x_{n-k+1}-1=0$. 
    This is enough to conclude that $x$ is the concatenation of slowly increasing sequences.
\end{proof}

Thanks to the above lemma, in a slow-climbing inversion sequence of an involution, the number of small ascents equals the number of nonzero elements. 
Another interesting consequence of Lemma~\ref{lem: relevant = slowly increasing} is the following corollary. 

\begin{corollary}\label{meet}
    Let $v,w\in \invol_n$ be two slow-climbing involutions of size $n$. Then $v\wedge w\in \invol_n$ is a slow-climbing involution as well.
\end{corollary}

\begin{proof}
	It is straightforward to check that taking the coordinate-wise minimum of two sequences that are concatenations of slowly increasing sequences results into a concatenation of slowly increasing sequences.  
\end{proof}

The above property can obviously be extended to any (finite) 
subset of slow-climbing involutions of the same size. Moreover, notice that Corollary \ref{meet} implies that the meet of slow-climbing involutions is still an involution. This fact will be crucial in the proof of Theorem \ref{moebius_computation}. Note that this property of meets is not true in general for arbitrary involutions, as $\invol_n$ is not a lattice for sufficiently large $n$. 

Given $w\in \invol_n$, 
let 
$M_w$
be the set of maximal slow-climbing 
involutions inside the principal order ideal $[\id ,w]$ in $\invol_n$. Clearly $M_w$ is an antichain, by definition. 

\begin{lemma}\label{nonzero_element}
    Let  $w\in \invol_n$ and suppose that $w\neq \id$. 
    Then there exists an index $i$ such that, for all $v\in M_w$, the $i$th component of $I(v)$ is nonzero. Equivalently, $\bigwedge M_w\neq \id$. 
\end{lemma}

\begin{proof}
    For an inversion sequence $y$, define a \emph{cluster} of $y$ to be any interval $[a,b] \subseteq [1,n]$ of indices, with the property that $y_{a+j} \ge j$ for all $j\in [0,b-a]$, and such that neither $[a-1,b]$ nor $[a,b+1]$ has this property. Equivalently, a cluster is an interval with the said property and whose width $b-a$ is maximal. It follows easily from this definition that any two clusters are incomparable with respect to containment. Moreover, the union of all clusters of $y$ is equal to $[1,n]$.

    Set $x:=I(w)=(x_1,\ldots ,x_n)$. 
  Without loss of generality, let us assume that $x_n > 0$ (otherwise, we can consider $(x_1,\ldots ,x_{n-1})$ and prove the statement for one size smaller). Let $C$ be the cluster of $x$ containing $n$. Note that there is a unique such cluster, by the maximality of clusters. Let $s$ be the number of elements in this cluster, and note that $s > 1$ because $x_n > 0$. 
    
    Suppose, for the purpose of obtaining a contradiction, that $n-1$ is the largest element of a cluster $D$ having $q$ elements. In fact, we must have $s \le q$ because clusters are incomparable with respect to containment. Moreover, $s-1 \le x_n < q$, by definition of the clusters $C$ and $D$. Since $x$ is the inversion sequence of an involution, Lemma~\ref{rec_char_inv} implies that the element in position $n-x_n$ of $x$ is equal to $0$. However, we have $n-x_n \in D$ but $n-x_n$ is not the smallest element of this cluster because $x_n < q$, and thus $x_{n-x_n}\neq 0$. This is a contradiction. Thus there is no such $D$, meaning that $C$ is the only cluster containing $n-1$. 
    
    Given $v \in [\id,w]$, assume that $v$ is slow-climbing and set $y := I(v) = (y_1, \ldots, y_n)$. If $y_n=0$, then the sequence obtained from $y$ by replacing $y_n$ with $y_{n-1}+1$ is slow-climbing, and the associated involution $\tilde{v}$ is such that $v < \tilde{v} \le w$, due to the fact that $n-1$ is included in no other cluster besides $C$. This implies that, for any $v \in M_w$, necessarily a maximal slow-climbing involution, we have $y_n > 0$. Therefore $\bigwedge M_w \neq \id$.
\end{proof}

The above proof gives a simple way to detect the value of $i$ in the statement of the lemma; namely, the said index $i$ is the position of the rightmost nonzero element of $x$.  We can now state and prove the main result of this section.

\begin{theorem}\label{moebius_computation}
	Fix an involution $w\in \invol_n$ 
 and let $\alpha$ be the number of nonzero entries in $\inv(w)$. Then 
	\[
	\mu (\id ,w)=
	\begin{cases}
	(-1)^\alpha &\textnormal{ if $w$ is slow-climbing, and} \\
	0 &\textnormal{ otherwise.}
	\end{cases}
	\]
\end{theorem}

\begin{proof}
	Set $x := \inv(w) = (x_1, \ldots, x_n)$. The proof is by induction on $|x|:=x_1 +\cdots +x_n$.
	
	If $|x|=0$, then $w=\id$ and $\mu (\id ,\id )=(-1)^0=1$.
	
	Now assume that $|x|\neq 0$ (equivalently, that $w\neq \id$), and suppose the result holds for involutions whose inversion sequences sum to smaller values. We start by supposing that $w$ is not slow-climbing. Let $M:= M_w$ denote the antichain of maximal slow-climbing involutions inside $[\id ,w]$. In particular, $w\notin M$. The usual recursion for M\"obius function can be written as 
	\begin{equation}\label{irrelevant}
	\mu (\id ,w)=-\sum_{v\in \downarrow M}\mu (\id ,v)-\sum_{\substack{\id\leq v<w\\ v\notin \downarrow M}}\!\!\mu (\id ,v),
	\end{equation}
	where $\downarrow\!S$ denotes the order ideal (downset) generated by a subset $S$ in $\invol_n$. Each summand of the second sum is equal to 0 (since, by definition of $M$, each $v\notin \downarrow\!M$ is not slow-climbing and we can then use the induction hypothesis), hence the second sum vanishes. To compute the first sum, first observe that obviously $\downarrow\!M=\bigcup_{u\in M} \downarrow\!u$, where we write $\downarrow\!u := \downarrow\!\{u\}$ to ease notation. Thus, if we suppose that $|M|=k$, 
    then we can use inclusion-exclusion to get
	\[
	\sum_{v\in \downarrow M}\mu (\id ,v)=\sum_{v\in \bigcup\limits_{u\in M} \downarrow u}\!\!\mu (\id ,v)=\sum_{h=1}^{k}(-1)^{h-1} \sum_{\substack{T\subseteq M\\ |T|=h}} \ \sum_{v\in \bigcap\limits_{u\in T}\downarrow u}\!\!\mu (\id ,v).  
	\]
	If $T=\{ u_1,\ldots ,u_h\}$ (recalling that each $u_i \in M$), then Corollary~\ref{meet} implies that 
 $$\bigcap_{u\in T}\downarrow\!u=\downarrow\!u_1 \cap \cdots \cap \downarrow\!u_h =\downarrow\!(u_1\wedge \cdots \wedge u_h ),$$
since $u_1\wedge \cdots \wedge u_h \in \invol_n$. Therefore each $\sum_{v\in \bigcap_{u\in T}\downarrow\!u}\mu (\id ,v)$ is equal to zero (because they are all sums of the values $\mu (\id ,v)$ for all $v$ in a non-trivial interval), unless $u_1 \wedge \cdots \wedge u_h =\id$. However, this never happens, thanks to Lemma~\ref{nonzero_element}. Thus the whole first sum in \eqref{irrelevant} is equal to zero, and the theorem holds when $w$ is not slow-climbing.

	Now suppose that $w$ is slow-climbing, and set $x=(x_1 ,\ldots ,x_n)=\vartheta_1 \vartheta_2\cdots \vartheta_k$, where (with a little abuse of notations) $\vartheta_i =(0,1,\ldots ,h_i)$. Observe that, with such notations, the number of nonzero elements of $x$ is $n-k$. Clearly we have that
	\[
	\mu (\id ,w)=-\hspace{-.4in}\sum_{\substack{\id\leq v<w\\ \textnormal{$v$ is not slow-climbing}}}\hspace{-.4in}\mu (\id ,v) \ - \ \hspace{-.3in}\sum_{\substack{\id\leq v<w\\ \textnormal{$v$ is slow-climbing}}}\hspace{-.3in}\mu (\id ,v),
	\]
	and each term of the first sum is equal to 0 by the induction hypothesis. So, in order to compute $\mu (\id ,w)$, we need to evaluate the contributions $\mu (\id ,v)$ for all slow-climbing involutions $v<w$. Again by the inductive hypothesis, such contributions are precisely $(-1)^l$, where $l$ is the number of nonzero elements of $z:= I(v)=(z_1,\ldots ,z_n )$. Therefore what is left to do is to count how many slow-climbing inversion sequences of involutions $z<x$ there are with $l$ nonzero elements. This can be done using a simple ``balls and walls'' argument. 
 Namely, we can encode $x$ using the alphabet $\{ \circ ,| \}$, by replacing each element with $\circ$ and inserting $|$ between any two consecutive slowly increasing sequences. In such a way, $x$ becomes a word $\psi_x$ of length $n+k-1$, having $n$ balls and $k-1$ walls. Any slow-climbing inversion sequence of an involution $z<x$ can be obtained from $x$ by splitting one or more slowly increasing sequences into shorter ones. Therefore $z$ is encoded by a word $\psi_z$ that is obtained from $\psi_x$ by inserting one or more walls between consecutive balls. Since there are $n+k-2$ internal spaces (between any two symbols) in $\psi_x$, and each of the $k-1$ walls makes two of them unavailable, we have a total of $n+k-2-2(k-1)=n-k$ admissible spaces to insert new walls. If we require $z$ to have $l$ nonzero elements, then this means that we want to insert $n-k-l$ new walls. Summing up, and observing that $0\leq l<n-k$, we thus have
	\[
	-\hspace{-.2in}\sum_{\substack{\id\leq v<w\\ \textnormal{$v$ is slow-climbing}}}\hspace{-.3in}\mu (\id ,v)=-\hspace{-.1in}\sum_{l=0}^{n-k-1}\binom{n-k}{n-k-l}(-1)^l=-((1-1)^{n-k}-(-1)^{n-k})=(-1)^{n-k},
	\]
	yielding the desired formula.       
\end{proof}

As we observed at the end of Section~\ref{sec:counting}, the boolean principal order ideals of the poset $\poset_n$ are those of the form $[\id ,w]$, where $I(w)$ contains only 0 and 1. In such a case, $w$ is an involution if and only if $I(w)$ does not contain consecutive $1$s (due to Proposition~\ref{rec_char_inv}). 
All such involutions $w$ are slow-climbing (this can be seen as a consequence of Lemma~\ref{lem: relevant = slowly increasing} for instance).
If $I(w)$ does not contain consecutive $1$s, then every consecutive subsequence of $I(w)$ has the same property, hence the interval $[(0,\ldots,0),I(w)]$ consists entirely of inversion sequences of involutions. As a consequence, the M\"obius function of the pair $(\id ,w)$ in $\poset_n$ coincides with the M\"obius function of $(\id, w)$ in $\invol_n$. This is coherent with our results, since Theorem~\ref{moebius_computation} implies that such a M\"obius function is given by the number of nonzero entries of $I(w)$; i.e., the number of $1$s of $I(w)$, which is indeed what we already knew from the end of Section~\ref{sec:counting}.

\section{Further directions}\label{open_problems}

The main purpose of the present paper was to present the basic combinatorial properties of the middle order on permutations. For this reason, we have focused on topics that are often studied about combinatorially defined posets. We have stated and proved several structural and combinatorial properties, and laid the combinatorial foundations of the middle order. On the other hand, we claim that such order may have applications in other fields, but we have left them to future investigations. Here we sketch two possible avenues of research that are particularly intriguing to us, with many other potentially interesting directions also available.

\subsection{A Heyting algebra structure on permutations}

A \emph{Heyting algebra} is a lattice $\mathcal{H}$ with minimum $\0$ and maximum $\1$, such that the ``relative pseudocomplement'' of $x$ with respect to $y$ exists for all $x, y\in \mathcal{H}$. By definition, the \emph{relative pseudocomplement} of $x$ with respect to $y$ is the element
$$
x \rightsquigarrow y \ :=\max \{ z\in \mathcal{H}\, |\, x\wedge z\leq y\} .
$$

In a Heyting algebra $\mathcal{H}$, two important notions are pseudocomplements and regular elements. The \emph{pseudocomplement} of $x$ is defined by $\sim\!x:=x\rightsquigarrow \0$. It can be shown
that $x\leq \sim \sim\!x$. The reverse inequality, 
however, does not hold in general. An element $x$ of $\mathcal{H}$ is \emph{regular} whenever $x=\sim \sim\!x$. The subposet of regular elements of a Heyting algebra forms a
boolean algebra. 
For the interested reader, some useful textbooks on lattice theory where Heyting algebras are discussed in some detail are \cite{BD,Bly}. 

It is well known that every finite distributive lattice has a canonical Heyting algebra structure. 
It may then be interesting to explore the Heyting algebra structure of $\poset_n$. The next theorem gives a characterization of the relative pseudocomplement and pseudocomplement operations, and of regular elements, in terms of parameters on permutations and their inversion sequences.

\begin{theorem}\label{thm:Heyting}
    Fix permutations $v,w\in S_n$, with inversion sequences $\inv (v)=(x_1,\ldots ,x_n)$ and $\inv (w)=(y_1,\ldots ,y_n)$. Then $\inv (v\rightsquigarrow w)=(z_1,\ldots ,z_n)$, where
    $$
    z_i = \begin{cases}i-1 & \text{if } x_i \leq y_i, \text{ and}\\ y_i & \text{if } x_i > y_i. \end{cases}
    $$
    As a consequence, $\sim\!v$ is the $132$- and $231$-avoiding permutation obtained from $v$ by listing its right-to-left minima in decreasing order, followed by all the remaining elements in increasing order.

    A permutation $v$ is regular if and only if $v$ avoids $132$ and $231$.
\end{theorem}

\begin{proof}
    The computation of $\inv (v\rightsquigarrow w)$ easily follows from \eqref{eq:meet_and_join2}. 
    As a consequence, we get that $\inv (\sim\!v)=(\overline{x}_1,\ldots ,\overline{x}_n)$, where
    $$
    \overline{x}_i =\begin{cases} i-1 & \text{if } x_i =0, \text{ and}\\ 0 & \text{if } x_i \neq 0. \end{cases}
    $$
    This implies that the one-line notation of $\sim\!v$ consists of a decreasing sequence followed by a (possibly empty) increasing sequence, whose decreasing prefix consists of those elements corresponding to the $0$s in the inversion sequence of $v$ (i.e., as we saw in the proof of Theorem~\ref{thm:formula for b(n,k)}, the right-to-left minima of $v$) and whose increasing suffix consists of all remaining elements. Since such permutations avoid 132 and 231, the result follows. Finally, $v$ is regular if and only if, for every $i$, either $x_i=0$ or $x_i =i-1$; that is, by the above argument, if and only if $v$ avoids 132 and 231. 
\end{proof}

For example, consider the permutations $v=361592784$ and $w=614928537$.
We have $I(v) = (0,0,2,0,2,4,1,1,4)$ and  $I(w) = (0,0,0,2,1,5,0,3,5)$.
Then 
$$I(v\rightsquigarrow w)=(0,1,0,3,1,5,0,7,8),$$
so that $v\rightsquigarrow w=986421537$. Moreover, $\sim\!v=421356789$, which is a regular element.  

Heyting algebras are especially relevant in mathematical logic.
Indeed, in the same way that boolean algebras are the algebraic counterpart of classical propositional logic, Heyting algebras are the algebraic counterpart of intuitionistic logic \cite{BdJ}. It is then natural to ask if our Heyting algebras of permutations have any interesting logical interpretation in the context of intuitionistic logic. 

\subsection{A poset structure on parking functions}

A \emph{parking function of size $n$} is a sequence $(p_1 ,p_2 ,\ldots ,p_n )$ of positive integers smaller than or equal to $n$ such that, denoting with $(q_1 ,q_2 ,\ldots ,q_n )$ its nondecreasing rearrangement, we have that $q_i \leq i$, for all $i$. Nondecreasing 
inversion sequences can be interpreted as nondecreasing parking functions (after adding 1 to each entry of the sequence). The name ``parking function'' originates from a parking model, where $n$ cars try to park into $n$ labeled spots in a one-way street, and the preference of car $i$ is the spot labeled $p_i$. If the preferred spot is found to be occupied, then the car proceeds forward and parks in the first available spot, if it exists, otherwise it leaves the street. A parking function is a sequence of preferences that allows all cars to park.
Since it is known that every rearrangement of a parking function is still a parking function, 
equivalently we have that the set $PF_n$ of parking functions of size $n$ is (essentially, i.e., up to the shift by $1$ of every entry) the set of all rearrangements of all inversion sequences. We can endow $PF_n$ with a partial order structure, again using the coordinate-wise order, thus obtaining a poset with minimum $(1,1,\ldots 1)$,  and whose set of maximal elements is the set of rearrangements of  $\{ 1,2,\ldots n\}$. If we introduce a maximum element $\top$, then the resulting poset $\mathcal{PF}_n$ appears to have a more complicated structure than $\poset_n$, and it may be interesting to explore its properties. We close with a simple result stating some basic facts about $\mathcal{PF}_n$.

\begin{theorem}
    The poset $\mathcal{PF}_n$ is a lattice, which is neither distributive nor modular (when $n\geq 3$).
\end{theorem}

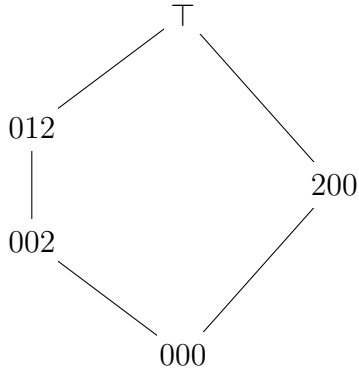
\begin{figure}[htbp]
\begin{center}
\begin{tikzpicture}
    \draw (0,0) coordinate (000);
    \draw (2,2.25) coordinate (200);
    \draw (-2,1.5) coordinate (002);
    \draw (-2,3) coordinate (012);
    \draw (0,4.5) coordinate (top);
    \draw (000) -- (200) -- (top) -- (012) -- (002) -- (000);
    \foreach \x in {000,200,002,012} {\fill[white] (\x)++(-.35,-.3) rectangle ++(.7,.6); \draw (\x) node {$\x$};}
    \fill[white] (top)++(-.25,-.25) rectangle ++(.5,.5);
    \draw (top) node {$\top$};
\end{tikzpicture}
\end{center}
\caption{A sublattice of $\mathcal{PF}_n$}\label{not_modular}
\end{figure}

\begin{proof}
    Given parking functions $p,q\in PF_n$, their coordinate-wise meet is still a parking function because its increasing rearrangement is less than or equal to the increasing rearrangement of both $p$ and $q$. Thus $\mathcal{PF}_n$ is a meet-semilattice and, since it has a maximum, it is also a lattice. To prove that it is not modular when $n\geq 3$ (and so, a fortiori, also not distributive), it is enough to find a sublattice isomorphic to the pentagon lattice. Figure~\ref{not_modular} shows such a sublattice when $n=3$, and the generalization for larger $n$ is immediate, for instance adding the appropriate number of final zeros to all elements.  
\end{proof}

\subsection*{Acknowledgments.} 

We are grateful to the organizers of the Dagstuhl Seminar \emph{Pattern Avoidance, Statistical Mechanics and Computational Complexity} and of the conference \emph{Permutation Patterns 2023} in Dijon, where this project started and developed, as well as to participants of these events who brought reference~\cite{AP} to our attention. 
We also thank Robert Johnson and Nick Ovenhouse, who pointed us to references \cite{JLL} and \cite{Nick} respectively, after a first version of our work was posted on the arXiv.

\end{document}